\documentstyle[12pt]{article}
\date{}

\begin{document}

\title{$H^{\infty}$ Performance of Interval Systems}
\author{Long Wang\thanks{Supported by National Natural Science Foundation of China (69925307), National Key Basic Research Special Fund (No. G1998020302) and the Research Fund for the Doctoral Program of Higher Education. }\\{\small Center for Systems and Control, Department of Mechanics and Engineering Science,}\\
{\small Peking University, Beijing 100871, CHINA}}
\maketitle

\begin{abstract}
In this paper, we study $H^{\infty}$ performance of interval
systems. We prove that, for an interval system, the maximal
$H^{\infty}$ norm of its sensitivity function is achieved at
twelve (out of sixteen) Kharitonov vertices.

\end{abstract}

\vspace{3mm}
\section{Introduction}

\hspace{1cm} Motivated by the seminal theorem of Kharitonov on robust stability of interval
polynomials\cite{Khar, Khar1}, a number of papers on robustness analysis of uncertain
systems have been published in the past few years\cite{Holl, Bart, Fu, wang1,
wang2, Barmish, Chap, wang3}. Kharitonov's theorem states that the Hurwitz stability of
the real (or complex) interval polynomial family can be guaranteed by the Hurwitz
stability of four (or eight) prescribed critical vertex polynomials in this family.
This result is significant since it reduces checking stability of infinitely many
polynomials to checking stability of finitely many polynomials, and the number of
critical vertex polynomials need to be checked is independent of the order of the
polynomial family. An important extension of Kharitonov's theorem is the edge theorem
discovered by Bartlett, Hollot and Huang\cite{Bart}. The edge theorem states that
the stability of a polytope of polynomials can be guaranteed by the stability of its
one-dimensional exposed edge polynomials. The significance of the edge theorem is
that it allows some (affine) dependency among polynomial coefficients, and applies to
more general stability regions, e.g., unit circle, left sector, shifted half plane,
hyperbola region, etc. When the dependency among polynomial coefficients is nonlinear,
however, Ackermann shows that checking a subset of a polynomial family generally can
not guarantee the stability of the entire family\cite{Ack1, Ack2}.

In this paper, we prove that, for an interval system, the maximal
$H_{\infty}$ norm of its sensitivity function is achieved at
twelve (out of sixteen) Kharitonov vertices. This result is useful
in robust performance analysis and $H_\infty$ control design for
dynamic systems under parametric perturbations.

\vspace{6mm}
\section{Main Results}

Denote the $m$-th, $n$-th $(m<n)$ order real interval polynomial families $K_g (s)$, $K_f (s)$ as

\begin{equation}
K_g(s)=\{g(s)|g(s)=\sum_{i=0}^mb_is^i, b_i\in [\underline{b_i},
\overline{b_i}], i=0,1, ......,m\},
\end{equation}

\begin{equation}
K_f(s)=\{f(s)|f(s)=\sum_{i=0}^na_is^i, a_i\in [\underline{a_i},
\overline{a_i}], i=0,1, ......,n\}.
\end{equation}

\noindent For any $f(s)\in K_f(s)$ , it can be expressed as

\begin{equation}
f(s)=\alpha _f(s^2)+s\beta _f(s^2),
\end{equation}

\noindent where

\begin{equation}
\alpha _f(s^2)=a_0+a_2s^2+a_4s^4+a_6s^6+......,
\end{equation}

\begin{equation}
\beta _f(s^2)=a_1+a_3s^2+a_5s^4+a_7s^6+.......
\end{equation}

\noindent Obviously, for any fixed $\omega \in R$, $\alpha _f(-\omega ^2)$ and $\omega \beta _f(-\omega^2)$ are the real and imaginary parts of $f(j\omega )\in C$ respectively.

For the interval polynomial family $K_f (s)$, define

\begin{equation}
\alpha _f^{(1)}(s^2)=\underline{a_0}+\overline{a_2}s^2+\underline{a_4}s^4+%
\overline{a_6}s^6+......,
\end{equation}

\begin{equation}
\alpha _f^{(2)}(s^2)=\overline{a_0}+\underline{a_2}s^2+\overline{a_4}s^4+%
\underline{a_6}s^6+......,
\end{equation}

\begin{equation}
\beta _f^{(1)}(s^2)=\underline{a_1}+\overline{a_3}s^2+\underline{a_5}s^4+%
\overline{a_7}s^6+......,
\end{equation}

\begin{equation}
\beta _f^{(2)}(s^2)=\overline{a_1}+\underline{a_3}s^2+\overline{a_5}s^4+%
\underline{a_7}s^6+......,
\end{equation}

\noindent and denote the four Kharitonov vertex polynomials of $K_f (s)$ as

\begin{equation}
f_{ij}(s)=\alpha _f^{(i)}(s^2)+s\beta _f^{(j)}(s^2), \,\,\,\,\,\,\,\,\,\,\,\, i,j=1,2
\end{equation}

For the interval polynomial family $K_g (s)$, the corresponding $\alpha _g^{(i)} (s), \beta _g^{(j)} (s)$ and $g_{ij}(s) \in K_g (s)$ can be defined analogously.

Denote by H the set of all Hurwitz stable polynomials (i.e.\, all of their roots lie within the open left half of the complex plane).

For the proper stable rational function $\frac{p(s)}{q(s)}$\,, the $H_{\infty}$ norm is defined as

\begin{equation}
  {\left \| \frac{p(s)}{q(s)} \right \|}_{\infty} = \sup \left \{ \left | \frac{p(j \omega )}{q(j \omega)} \right | \; \mid \omega \in (- \infty \;,\;+ \infty ) \right \}
\end{equation}

\noindent The proper complex rational function $\frac{p(s)}{q(s)}$ is said to be strictly positive real, if 1) $q(s) \in H$\,;\, and 2) for any $\omega \in R \;,\; \Re \;\frac{p(j \omega)}{q(j \omega)}>0$.

Denote by SPR the set of all strictly positive real rational functions.

\vspace{0.4cm}

\noindent {\bf Lemma 2.1}\cite{das}

For any fixed $\omega \in R$, $f(s)\in K_f(s)$, we have

\begin{equation}
\alpha _f^{(1)}(-\omega ^2)\leq \alpha _f(-\omega ^2)\leq \alpha
_f^{(2)}(-\omega ^2),
\end{equation}

\begin{equation}
\beta _f^{(1)}(-\omega ^2)\leq \beta _f(-\omega ^2)\leq \beta
_f^{(2)}(-\omega ^2).
\end{equation}

\vspace{0.4cm}

\noindent {\bf Lemma 2.2}\cite{Ack1} (Zero Exclusion Principle)

For the $n$-th order polynomial family

\begin{equation}
f(s,T)=: \{f(s,t)|t\in T\},
\end{equation}

\noindent where $T$ is a bounded connected closed set, and the coefficients of $f(s,t)$ are continuous functions of $t$, then $f(s,T)\in H$ if and only if

1)\,\,\,\,\,there exists $t^{*}\in T,$  such that $f(s,t^{*})\in H$;

2)\,\,\,\,\,$0\notin f(j\omega ,T),$ $\forall \omega \in R$.

Consider the strictly proper open-loop transfer function

\begin{equation}
P=\frac{g(s)}{f(s)}
\end{equation}

\noindent and suppose the closed-loop system is stable under negative unity feedback. Denote its sensitivity function as

\begin{equation}
S=\frac{1}{1+P}=\frac{f(s)}{f(s)+g(s)}
\end{equation}

\noindent Apparently, we have

\begin{equation}
||S||_{\infty} \geq 1
\end{equation}

For notational simplicity, define

\begin{equation}
J_{i_1 j_1 i_2 j_2} (s) = g_{i_1 j_1}(s) + (1 + \delta e^{j \theta}) f_{i_2 j_2}(s), \,\,\,\,\,\,\,\,\,\, \delta \in (0, 1), \,\,\, i_1, j_1, i_2, j_2 = 1,2, \,\,\, \theta \in [-\pi, \pi].
\end{equation}

\vspace{0.4cm}

\noindent {\bf Lemma 2.3}

Suppose $g(s) + f(s) \in H$. Then, for any $\gamma > 1$, we have

\begin{equation}
||S||_{\infty}<\gamma \Longleftrightarrow g(s) + (1 + \frac{1}{\gamma} e^{j \theta}) f(s) \in H, \,\,\,\,\, \forall \theta \in [-\pi, \pi].
\end{equation}

\vspace{0.6cm}

\noindent Proof: Necessity: Since $g(s) + f(s) \in H$ and $|| \frac{ \frac{1}{\gamma} f(s)}{f(s) + g(s)} ||_{\infty} < 1$, by Rouche's Theorem, we know that

\begin{equation}
[g(s) +f(s)] + \frac{1}{\gamma} e^{j \theta} f(s) \in H, \,\,\,\,\, \forall \theta \in [-\pi, \pi]
\end{equation}

Sufficiency: Now suppose on the contrary that $||S||_{\infty} \geq \gamma$, namely, $|| \frac{ \frac{1}{\gamma} f(s)}{f(s) + g(s)} ||_{\infty} \geq 1$. Since
$| \frac{ \frac{1}{\gamma} f(s)}{f(s) + g(s)} |_{s=j \omega} |$ is a contiunous function of $\omega$, and since

\begin{equation}
\lim_{\omega \rightarrow \infty} | \frac{ \frac{1}{\gamma} f(s)}{f(s) + g(s)} |_{s=j \omega} | = \frac{1}{\gamma} < 1
\end{equation}

\noindent there must exist $\omega_0$ such that

\begin{equation}
| \frac{ \frac{1}{\gamma} f(s)}{f(s) + g(s)} |_{s=j \omega_0} | = 1
\end{equation}

\noindent Therefore, there exists $\theta_0 \in [-\pi, \pi]$ such that

\begin{equation}
\{ g(s) +f(s) + \frac{1}{\gamma} e^{j \theta_0} f(s) \} |_{s = j \omega_0} = 0
\end{equation}

\noindent which contradicts the original hypothesis. This completes the proof.

\vspace{0.4cm}

\noindent {\bf Lemma 2.4}

For any $\delta \in (0, 1), \theta \in [-\pi, \pi]$, we have

\begin{equation}
W(s)=: \{g(s) + (1 + \delta e^{j \theta}) f(s) | g(s)\in K_g(s), f(s)\in K_f(s) \} \subset H \Longleftrightarrow
\end{equation}

\begin{equation}
J_{1111}, J_{1212}, J_{2222}, J_{2121}, J_{1112}, J_{1222}, J_{2221}, J_{2111}, J_{1211}, J_{2212}, J_{2122}, J_{1121} \in H
\end{equation}

\vspace{0.6cm}

\noindent Proof: Necessity is obvious. To prove sufficiency, note that $W(s)$ is a set of polynomials with complex coefficients, and with constant order $n$.
By Lemma 2.2, it suffices to show that

\begin{equation}
0 \not\in W(j \omega), \,\,\,\,\, \forall \omega \in R
\end{equation}

\noindent Since $0 \not\in W(j \omega_{\infty})$ for sufficiently large $\omega_{\infty}$, we only need to show that

\begin{equation}
0 \not\in \partial W(j \omega), \,\,\,\,\, \forall \omega \in R
\end{equation}

\noindent where $\partial W(j \omega)$ stands for the boundary of $W(j \omega)$ in the complex plane.

To construct $\partial W(j \omega)$, note that $\arg (1+ \delta e^{j \theta}) \in (- \frac{\pi}{2}, \frac{\pi}{2})$. Suppose now
$\omega \geq 0$ and $\arg (1+ \delta e^{j \theta}) \in [0, \frac{\pi}{2})$. Then by Lemma 2.1, we know that
$K_g (j \omega), K_f (j \omega)$ are rectangles with edges parallel to the coordinate axes. The four vertices of
$K_g (j \omega)$ are $g_{11} (j \omega), g_{12} (j \omega), g_{21} (j \omega), g_{22} (j \omega)$, respectively;
and the four vertices of $K_f (j \omega)$ are $f_{11} (j \omega),$ $f_{12} (j \omega),$ $f_{21} (j \omega), f_{22} (j \omega)$, respectively.
$(1+ \delta e^{j \theta}) K_f (j \omega)$ is generated by rotating $K_f (j \omega)$ by
$\arg (1+ \delta e^{j \theta})$ counterclockwisely, and then scaling by $|1+ \delta e^{j \theta}|$. Thus,
$W(j \omega) = K_g (j \omega) + (1+ \delta e^{j \theta}) K_f (j \omega)$ is a convex polygon with eight edges.
These edges are parallel to either the edges of $K_g (j \omega)$ or the edges of $(1+ \delta e^{j \theta}) K_f (j \omega)$.
Therefore, their orientations are fixed (independent of $\omega$). The eight vertices of $W(j \omega)$ are (clockwisely) $J_{1111} (j \omega)$,
$J_{1112} (j \omega)$, $J_{1212} (j \omega)$, $J_{1222} (j \omega)$, $J_{2222} (j \omega)$, $J_{2221} (j \omega)$, $J_{2121} (j \omega)$, $J_{2111} (j \omega)$, respectively.

Now suppose on the contrary that there exists $\omega_0 \geq 0$ such that

\begin{equation}
0 \in \partial W(j \omega_0)
\end{equation}

\noindent Without loss of generality, suppose

\begin{equation}
0 \in \{ \lambda J_{1111} (j \omega_0) + (1- \lambda) J_{1112} (j \omega_0) | \lambda \in [0, 1] \}
\end{equation}

\noindent Namely, there exists $\lambda_0 \in (0, 1)$ such that

\begin{equation}
\lambda_0 J_{1111} (j \omega_0) + (1- \lambda_0) J_{1112} (j \omega_0) = 0
\end{equation}

\noindent Since $J_{1111} (s), J_{1112} (s) \in H$, we have

\begin{equation}
\frac{d}{d \omega} \arg J_{1111} (j \omega) > 0, \,\,\,\,\,\,\,\, \frac{d}{d \omega} \arg J_{1112} (j \omega) > 0
\end{equation}

\noindent Thus\cite{Ran}

\begin{equation}
\hspace{-5cm} \frac{d}{d \omega} \arg [J_{1112} (j \omega) - J_{1111} (j \omega)] |_{\omega = \omega_0} =
\end{equation}

\begin{equation}
(1- \lambda_0) \frac{d}{d \omega} \arg J_{1111} (j \omega) |_{\omega = \omega_0} + \lambda_0 \frac{d}{d \omega} \arg J_{1112} (j \omega) |_{\omega = \omega_0} > 0
\end{equation}

\noindent This contradicts the fact that the edges of $W(j \omega)$ have fixed orientations. Thus

\begin{equation}
0 \not\in \partial W(j \omega)
\end{equation}

\noindent Suppose now
$\omega \leq 0$ and $\arg (1+ \delta e^{j \theta}) \in (- \frac{\pi}{2}, 0]$. Then
$K_g (j \omega), (1+ \delta e^{j \theta}) K_f (j \omega)$ are the mirror images (with respect to the real axis) of the corresponding sets
in the case of $\omega \geq 0$ and $\arg (1+ \delta e^{j \theta}) \in [0, \frac{\pi}{2})$. Therefore, following an identical line of arguments, we have

\begin{equation}
0 \not\in \partial W(j \omega)
\end{equation}

\noindent The cases when $\omega \geq 0$ and $\arg (1+ \delta e^{j \theta}) \in (- \frac{\pi}{2}, 0]$ and when
$\omega \leq 0$ and $\arg (1+ \delta e^{j \theta}) \in [0, \frac{\pi}{2})$ are also symmetric with respect to the real axis.
Hence, we only need to consider the former case. In this case, $K_g (j \omega), K_f (j \omega)$ are rectangles with edges parallel to the coordinate axes.
$(1+ \delta e^{j \theta}) K_f (j \omega)$ is generated by rotating $K_f (j \omega)$ by
$| \arg (1+ \delta e^{j \theta}) |$ clockwisely, and then scaling by $|1+ \delta e^{j \theta}|$. Thus,
$W(j \omega) = K_g (j \omega) + (1+ \delta e^{j \theta}) K_f (j \omega)$ is a convex polygon with eight edges.
These edges are parallel to either the edges of $K_g (j \omega)$ or the edges of $(1+ \delta e^{j \theta}) K_f (j \omega)$.
Therefore, their orientations are fixed (independent of $\omega$). The eight vertices of $W(j \omega)$ are (clockwisely) $J_{1111} (j \omega)$,
$J_{1211} (j \omega)$, $J_{1212} (j \omega)$, $J_{2212} (j \omega)$, $J_{2222} (j \omega)$, $J_{2122} (j \omega)$, $J_{2121} (j \omega)$, $J_{1121} (j \omega)$, respectively.
Thus, following a similar argument, we have

\begin{equation}
0 \not\in \partial W(j \omega)
\end{equation}

\noindent This completes the proof.

The following theorem shows that, for an interval system, the maximal $H_{\infty}$ norm of its
sensitivity function is achieved at twelve (out of sixteen) Kharitonov vertices.

\vspace{0.3cm}

\noindent {\bf Theorem 2.1}

Suppose $g_{ij} (s) + f_{ij} (s) \in H, \,\, i,j=1,2$. Then

\begin{equation}
\hspace{-5cm} \max \{ || \frac{f(s)}{f(s) + g(s)} ||_{\infty} | g(s)\in K_g(s), f(s)\in K_f(s) \} =
\end{equation}

\begin{equation}
\hspace{-4.5cm} \max \{ || \frac{f_{i_2 j_2} (s)}{f_{i_2 j_2} (s) + g_{i_1 j_1} (s)} ||_{\infty} | (i_1 j_1 i_2 j_2) = (1111), (1212),
\end{equation}

\begin{equation}
(2222), (2121), (1112), (1222), (2221), (2111), (1211), (2212), (2122), (1121) \}
\end{equation}

\vspace{0.6cm}

\noindent Proof: Since $g_{ij} (s) + f_{ij} (s) \in H, \,\, i,j=1,2$, by Kharitonov's Theorem\cite{Khar}, we know that
$K_g(s) + K_f(s) \subset H$. Let

\begin{equation}
\hspace{-5cm} \gamma_1 = \max \{ || \frac{f(s)}{f(s) + g(s)} ||_{\infty} | g(s)\in K_g(s), f(s)\in K_f(s) \}
\end{equation}

\begin{equation}
\hspace{-4cm} \gamma_2 = \max \{ || \frac{f_{i_2 j_2} (s)}{f_{i_2 j_2} (s) + g_{i_1 j_1} (s)} ||_{\infty} | (i_1 j_1 i_2 j_2) = (1111), (1212),
\end{equation}

\begin{equation}
(2222), (2121), (1112), (1222), (2221), (2111), (1211), (2212), (2122), (1121) \}
\end{equation}

\noindent Then apparently

\begin{equation}
\gamma_1 \geq \gamma_2 \geq 1
\end{equation}

Now suppose $\gamma_1 \neq \gamma_2$, namely, $\gamma_1 > \gamma_2$. Then there exists
$\gamma_0$ such that $\gamma_1 > \gamma_0 > \gamma_2$. Thus, for any
$(i_1 j_1 i_2 j_2) \in \{(1111), (1212), (2222), (2121), (1112), (1222), (2221),$ $(2111), (1211), (2212), (2122), (1121) \}$, we have

\begin{equation}
|| \frac{f_{i_2 j_2} (s)}{f_{i_2 j_2} (s) + g_{i_1 j_1} (s)} ||_{\infty} < \gamma_0
\end{equation}

\noindent Hence, by Lemma 2.3, we have

\begin{equation}
g_{i_1 j_1} (s) + (1 + \frac{1}{\gamma_0} e^{j \theta}) f_{i_2 j_2} (s) \in H, \,\,\,\,\, \forall \theta \in [-\pi, \pi]
\end{equation}

\noindent By Lemma 2.4, we know that

\begin{equation}
\{ g(s) + (1 + \frac{1}{\gamma_0} e^{j \theta}) f(s) | g(s)\in K_g(s), f(s)\in K_f(s) \} \subset H, \,\,\,\,\, \forall \theta \in [-\pi, \pi]
\end{equation}

\noindent Therefore, by Lemma 2.3, for any $g(s)\in K_g(s), f(s)\in K_f(s)$, we have

\begin{equation}
|| \frac{f(s)}{f(s) + g(s)} ||_{\infty} < \gamma_0
\end{equation}

\noindent Namely

\begin{equation}
\max \{ || \frac{f(s)}{f(s) + g(s)} ||_{\infty} | g(s)\in K_g(s), f(s)\in K_f(s) \} < \gamma_0
\end{equation}

\noindent That is, $\gamma_1 < \gamma_0$, which contradicts $\gamma_1 > \gamma_0 > \gamma_2$. This completes the proof.

\vspace{6mm}
\section{Conclusions}

We have proved that, for an interval system, the maximal
$H^{\infty}$ norm of its sensitivity function is achieved at
twelve (out of sixteen) Kharitonov vertices. This result is useful
in robust performance analysis and $H_\infty$ control design for
dynamic systems under parametric perturbations.

\end{document}